\documentclass{article}

\usepackage{arxiv}

\usepackage[utf8]{inputenc} 
\usepackage[T1]{fontenc}    
\usepackage{hyperref}       
\usepackage{url}            
\usepackage{booktabs}       
\usepackage{amsfonts}       
\usepackage{nicefrac}       
\usepackage{microtype}      
\usepackage{lipsum}		
\usepackage{graphicx}
\usepackage{doi}
\usepackage{amsthm}
\usepackage[square,sort,comma,numbers]{natbib}

\usepackage{amssymb,amsthm,amsmath}
\usepackage{xcolor,paralist,hyperref,titlesec,fancyhdr,etoolbox}
\usepackage{algorithm}
\usepackage{algpseudocode}

\title{Penalty-Free SDDP: Feasibility Cuts for Robust Multi-Stage Stochastic Optimization in Energy Planning}

\author{
\hspace{-2.5em}Guilherme~Freitas\\
\hspace{-2.5em}	PSR\\
\hspace{-2.5em}	Rio de Janeiro, Brazil\\
\hspace{-2.5em}	\texttt{guilhermefreitas@psr-inc.com} \\
\And
\hspace{-2.5em}Luiz~Carlos~da~Costa~Junior \\
\hspace{-2.5em}	PSR\\
\hspace{-2.5em}	Rio de Janeiro, Brazil\\
\hspace{-2.5em}	\texttt{luizcarlos@psr-inc.com} \\
\And
\hspace{-2.5em}Tiago~Andrade \\
\hspace{-2.5em}	PSR\\
\hspace{-2.5em}	Rio de Janeiro, Brazil\\
\hspace{-2.5em}	\texttt{tiago.andrade@psr-inc.com} \\
\And
Alexandre~Street \\
	Department of Electrical Engineering\\
	PUC-Rio\\
	Rio de Janeiro, Brazil\\
	\texttt{street@puc-rio.br} \\
}

\date{}

\renewcommand{\headeright}{}
\renewcommand{\shorttitle}{Penalty-Free SDDP: Feasibility Cuts for Robust Multi-Stage Stochastic Optimization in Energy Planning}

\hypersetup{
pdftitle={Penalty-Free SDDP: Feasibility Cuts for Robust Multi-Stage Stochastic Optimization in Energy Planning},
pdfsubject={},
pdfauthor={Guilherme Freitas, Luiz Carlos da Costa Junior, Tiago Andrade, Alexandre Street},
pdfkeywords={Penalty-Free SDDP, feasibility cuts, multi-stage stochastic optimization, Benders decomposition, hydrothermal scheduling, energy planning.},
}

\begin{document}
\maketitle

\begin{abstract}
	Multi-stage decision problems under uncertainty can be efficiently solved with the Stochastic Dual Dynamic Programming (SDDP) algorithm. However, traditional implementations require all stage problems to be feasible. Feasibility is usually enforced by adding slack variables and penalizing them in the objective function, a process that depends on case-specific calibration and often distorts the economic interpretation of results. This paper proposes the \textit{Penalty-Free SDDP}, an extension that introduces a \textit{Future Feasibility Function} alongside the traditional Future Cost Function. The new recursion handles infeasibilities automatically, distinguishing between temporary and \textit{truly infeasible} cases, and propagates feasibility information across stages through dedicated feasibility cuts. The approach was validated in a large-scale deterministic case inspired by the Brazilian hydrothermal system, achieving equivalent feasibility to the benchmark solution while eliminating miscalibrated artificial penalties. Results confirm its robustness and practicality as a foundation for future stochastic, multi-stage applications.
\end{abstract}


\section{Introduction} 
\label{section}

Multistage decision problems under uncertainty arise in several domains, including energy systems, finance, logistics, and water management. These problems require sequential decisions that depend on both current information and uncertain future outcomes. Solving them in their extensive form is computationally intractable because of the exponential growth of the scenario tree. To overcome this limitation, \textit{Stochastic Dual Dynamic Programming} (SDDP)~\cite{SDDP1991} is the state-of-the-art algorithm, relying on decomposition to approximate the cost-to-go functions, which represents a decision policy that guide optimal decisions over time.

A general assumption in traditional SDDP is that of \textit{relatively complete recourse}, which guarantees that every realization of uncertainty leads to a feasible problem. In practice, this assumption is rarely satisfied and feasibility is typically enforced by introducing slack variables penalized in the objective function. Although this artifice allows the modified problem to be numerically solved, it mixes real economic costs with artificial penalties, which lack clear interpretation and must be subjectively calibrated. The calibration process itself is not trivial, especially when violations involve quantities of different natures such as hydraulic limits, energy balance, or environmental targets. As a consequence, the modified problem can produce distorted price signals and misleading economic interpretations, ultimately leading to suboptimal planning decisions.

The situation is particularly problematic when constraint violations are associated with low-probability events. Their impact becomes diluted in the expected cost, making it difficult to guarantee constraint satisfaction unless penalties are increased to unrealistic values or risk-averse formulations are introduced, both of which further bias the results. In the former case, where penalties are artificially inflated to enforce compliance, convergence tends to deteriorate, and numerical instability may arise.

As modern planning problems grow in complexity, driven by renewable integration, regulatory requirements, and sustainability goals, there is a growing need for methodologies that can handle feasibility in a rational and economically consistent way.

This paper proposes a general extension to SDDP that eliminates the need for artificial penalties by introducing \textit{Benders feasibility cuts}. The proposed \textit{Penalty-Free SDDP} expands the classical algorithm with an additional recursive function, the \textit{Future Feasibility Function (FFF)}, which quantifies the value of preserving flexibility to meet future constraints. The method maintains the decomposition structure and convergence properties of SDDP, while ensuring model robustness and interpretability. Although the proposed framework is general and applicable to any multi-stage stochastic optimization problem, its performance is illustrated here through an application to long-term energy planning under operational constraints.

\section{Background and Related Work}

The foundation of the proposed approach lies in the decomposition principles introduced by Benders~\cite{Benders}, later extended to stochastic optimization by Van Slyke and Wets~\cite{LShaped}. In this framework, a large-scale problem is decomposed into a master problem and a set of subproblems. Information from the subproblems is then summarized by linear inequalities, known as \textit{Benders cuts}, that iteratively refine the master problem. Two types of cuts are typically generated: \textit{optimality cuts}, when the subproblem is feasible, and \textit{feasibility cuts}, when it becomes infeasible after the master’s decision.

Building upon this idea, Pereira and Pinto~\cite{SDDP1991} proposed the Stochastic Dual Dynamic Programming (SDDP) algorithm as a multi-stage extension of the L-shaped method. SDDP combines Benders-type decomposition with Monte Carlo sampling to approximate the stochastic process, recursively constructing a piecewise-linear representation of the expected future cost at each stage through a set of linear cuts, known as the \textit{Future Cost Function (FCF)}. This combination of sampling and decomposition enables the solution of complex multi-stage stochastic problems while avoiding the exponential growth of the scenario tree that would result from an explicit formulation.

Several works have investigated alternatives to the penalty-based approach, including the generation of feasibility cuts within the SDDP framework~\cite{Santos2009,Yuen,Guigues2014,HELSETH2023,HELSETH2024}. These studies show that feasibility cuts can accelerate convergence and improve numerical stability in problems that are ultimately feasible. However, most formulations rely on simplified or problem-specific assumptions, do not address the continuation of the algorithm in truly infeasible cases, or mix feasibility cuts with risk-based formulations that alter the original objective function. As a result, a general penalty-free extension of SDDP that maintains its decomposition structure and robustness has not yet been fully established~\cite{GUIGUES2023752,Diniz2020,SDDPReview}.

The present work addresses this gap by formulating a \textit{Penalty-Free SDDP} algorithm that explicitly integrates feasibility cuts into the recursive process while preserving the economic meaning of the cost function. This approach ensures feasibility handling without the need for artificial penalties and provides a consistent framework applicable to a broad class of multi-stage stochastic optimization problems.

\section{Penalty-Free SDDP Framework}

\subsection{Classical SDDP with penalized slacks}

In the classical penalty-based formulation, feasibility is ensured at each stage by introducing nonnegative slack variables that are penalized in the objective function. Let $x_t$ be the stage-$t$ decision, $s_t \ge 0$ the slack vector associated with potentially violated constraints, and let $w_t \ge 0$ denote the penalty weights. Using a generic stagewise functional $\mathcal{R}_t[\cdot]$ (e.g. expectation or risk measure), the stage problem can be formulated as

\begin{equation}
\label{eq:sddp_classic}
\begin{aligned}
\mathrm{FCF}_t(x_{t-1},&\xi_{t-1}) = \\
\min_{x_t,s_t} \;\; & c_t^\top x_t + w_t^\top s_t + \mathcal{R}_t\!\bigl[\widetilde{\mathrm{FCF}}_{t+1}(x_t,\xi_t)\bigr] \\
\text{s.t.} \;\; & A_t x_t + s_t \geq b_t - D_t x_{t-1},\\
& x_t \ge 0,\;\\ & s_t \ge 0.
\end{aligned}
\end{equation}

Here, $\widetilde{\mathrm{FCF}}_{t+1}$ is the piecewise-linear approximation of the Future Cost Function constructed from Benders optimality cuts. The penalty vector $w_t$ is chosen ``large enough'' to discourage violations, but its calibration is problem-dependent and lacks direct economic meaning. As discussed in the introduction, this artifice allows the modified problem to be solved at every iteration while blending economic costs with artificial penalties, which complicates interpretation and may affect numerical behavior, especially when low-probability violations are present.

\subsection{Penalty-Free SDDP formulation}

The proposed formulation extends the traditional SDDP recursion by splitting the solution process into two coordinated stages: a \textit{feasibility problem}, responsible for identifying and quantifying violations of the system constraints, and an \textit{optimality problem}, which computes the economic decision based on the feasible domain defined by the previous step. This separation allows the algorithm to treat feasibility and optimality independently, removing artificial penalties from the cost function while maintaining the decomposition structure of SDDP.

\vspace{2mm}
\noindent\textbf{Step 1: Feasibility problem (strict, all scenarios).}
At stage $t$, feasibility is enforced \emph{for all} realizations $\xi_t \in \Xi_t$. 
Nonnegative slacks $s_t$ absorb constraint violations and are weighted by $w_t \!>\! 0$ to express preferences when true infeasibilities occur. 
Unlike the penalized formulation of classical SDDP, these weights do not establish a financial trade-off between economic costs and constraint violations; instead, they serve only to prioritize among different types of violations (e.g., hydraulic, balance, or environmental) in situations where full feasibility cannot be achieved.

The feasibility recursion does not rely on probabilistic aggregation (such as expectation or CVaR). 
Instead, it follows a worst-case principle, equivalent to a conservative risk operator that enforces feasibility across all realizations. 
Algorithmically, this principle is implemented by retaining the \emph{union of feasibility cuts} encountered across sampled scenarios, guaranteeing that the solution remains feasible under every possible realization. 
A compact robust formulation is:

\begin{equation}
\label{eq:feasibility_problem_strict}
\begin{aligned}
\mathrm{FFF}_t(x_{t-1},&\xi_{t-1}) = \\
\min_{x_t,\, s_t,\, \beta_t} \quad & w_t^\top s_t + \beta_t \\[1mm]
\text{s.t.}\quad
    & A_t x_t + s_t \ge b_t - D_t x_{t-1}, \\
    & \beta_t \ge \widetilde{\mathrm{FFF}}_{t+1}(x_t,\xi_t), \quad \forall \xi_t \in \Xi_t, \\
    & x_t \ge 0, \\
    & s_t \ge 0.
\end{aligned}
\end{equation}

This problem computes the minimum infeasibility required to satisfy all constraints for a given state $x_{t-1}$, taking into account both immediate and future constraint violations. 
Its objective combines two terms: (i) the current infeasibility, measured by the weighted slacks $w_t^\top s_t$, and (ii) the recursive infeasibility $\beta_t$, which represents the value of the Future Feasibility Function (FFF) at the next stage. 
The solution provides two key quantities: the minimal slack vector $s_t^*$ and the cumulative infeasibility $\beta_t^*$. 
These values quantify the degree of constraint violation at the current and future stages, respectively, and are used to restrict the feasible region in the subsequent optimality problem. 
In particular, $s_t^*$ defines admissible bounds for local constraint relaxations, while $\beta_t^*$ limits the propagation of infeasibility in the recursive structure.

Each time a violation is observed in any scenario, a new feasibility cut is generated and retained, progressively constructing a piecewise-linear outer approximation of the feasible domain. 
This “union-of-cuts” approach ensures that both $s_t^*$ and $\beta_t^*$ capture all relevant infeasibilities, maintaining strict feasibility across forward and backward passes while preserving the stagewise decomposition structure of the algorithm.

\vspace{2mm}
\noindent\textbf{Step 2: Optimality problem with feasibility cuts.}
The second step solves the cost-minimization subproblem, where the previously computed infeasibility bounds $s_t^*$ and $\beta_t^*$ are enforced as upper limits for the slack and cumulative infeasibility variables, respectively. 
This formulation incorporates all feasibility information into the economic optimization while keeping the objective function free of artificial penalties.

Formally, the stage problem at time $t$ and realization $\xi_t$ is expressed as:
\begin{equation}
\label{eq:optimality_problem_with_FFF}
\begin{aligned}
{\mathrm{FCF}}_t(x_{t-1}, \xi&_{t-1}) = \\
\min_{x_t,\, s_t,\, \beta_t} \;
    & c_t^\top x_t + \mathcal{R}_t\!\big[ \widetilde{\mathrm{FCF}}_{t+1}(x_t, \xi_t) \big] \\
\text{s.t.} \;
    & A_t x_t + s_t \geq b_t - D_t x_{t-1}, \\
    & x_t \ge 0, \\
    & 0 \le s_t \le s_t^*, \\
    & \beta_t \le \beta_t^*, \\
    & \beta_t \ge \widetilde{\mathrm{FFF}}_{t+1}(x_t,\xi_t), \quad \forall \xi_t \in \Xi_t.
\end{aligned}
\end{equation}

When $\beta_t^* = 0$, the current state is in the approximated FFF domain. Then, the cumulative feasibility slack $\beta_t$ becomes inactive, and the problem coincides with the standard SDDP formulation with feasibility cuts. Otherwise, when $\beta_t^* > 0$, it propagates the residual infeasibility backward, allowing the model to adjust the current state variables to anticipate constraint violations. In some visited states, the problem may still be infeasible — a situation that arises when maintaining feasibility requires multi-stage anticipation or when a true structural infeasibility is present in the data but has not yet been fully detected.

The constraint $\beta_t \ge \widetilde{\mathrm{FFF}}_{t+1}(x_t,\xi_t)$ enforces consistency with the recursive feasibility representation, ensuring that any remaining violation at future stages is recognized at the current one. Through this coupling, the algorithm becomes \textit{feasibility-aware}: state decisions such as storage levels, thermal fuel contracts, emission budgets or other intertemporal controls are guided not only by expected economic trade-offs (captured by the FCF) but also by the accumulated infeasibility information (captured by the FFF).

The dual multipliers obtained from this problem generate \textit{optimality cuts} that refine the Future Cost Function, while the feasibility cuts remain active to define the feasible region. As the backward iterations progress, both sets of cuts are updated and expanded, progressively improving the joint representation of cost and feasibility across stages.

\section{Algorithmic implementation}

The Penalty-Free SDDP algorithm extends the standard forward–backward procedure by explicitly solving a feasibility problem before the optimality problem at each stage, obtaining the slack bounds $s_t^\ast$ and $\beta^*$ and enforcing them in the cost minimization. Convergence reflects the stabilization of both feasibility and cost approximations.

\subsection{Forward pass}

In the forward pass, a set of scenario paths $\Omega^{\mathrm{fwd}}$ is sampled. For each path $\omega \in \Omega^{\mathrm{fwd}}$ and for $t=1,\dots,T$, we visit the corresponding conditional state and solve, in this order:
\begin{enumerate}
    \item the \emph{feasibility problem}~\eqref{eq:feasibility_problem_strict} at $(x_{t-1}(\omega),\xi_{t-1}(\omega))$, obtaining the minimum slack vector $s_t^\ast(\omega)$ and the cumulative infeasibility $\beta_t^\ast(\omega)$;\vspace{-2pt}
    \item the \emph{optimality problem with feasibility cuts}~\eqref{eq:optimality_problem_with_FFF} at the same state, now enforcing $s_t \le s_t^\ast(\omega)$ and $\beta_t \le \beta_t^\ast(\omega)$.
\end{enumerate}
The forward pass therefore \emph{directly obtains} the violations (via $s_t^\ast$ and $\beta_t^\ast$) and \emph{visits new states} $x_t(\omega)$ in the conditional scenario tree, building the set of trial states $\mathcal{S}_t := \mathcal{S}_t \cup \{x_t(\omega)\}$ for each stage $t$.

\subsection{Backward pass}

In the backward pass, the algorithm proceeds from $t=T$ down to $t=1$. For each stage $t$ and for every trial state $x_t \in \mathcal{S}_t$ generated in the forward pass, we solve:
\begin{itemize}
    \item the feasibility problem~\eqref{eq:feasibility_problem_strict} (at $x_{t}$ and across the sampled realizations at $t{+}1$) obtaining the minimum slack vector $s_t^\ast$ and the cumulative infeasibility $\beta_t^\ast$; and, from its dual multipliers, \emph{generate feasibility cuts} that update $\widetilde{\mathrm{FFF}}_{t}$;
    \item the optimality problem~\eqref{eq:optimality_problem_with_FFF} (with $s_t\le s_t^\ast$ and $\beta \le \beta_t^\ast$ inherited from the feasibility problem above) and, from its dual multipliers, \emph{generate optimality cuts} that update $\widetilde{\mathrm{FCF}}_{t}$.
\end{itemize}
Both sets of cuts are stored, refining the piecewise-linear approximations of the FCF and FFF at each stage. When no new feasibility cuts are produced in the backward pass, the feasible region is considered stable at the trial states visited.

\subsection{Convergence criteria}

Convergence is declared when the following conditions hold simultaneously:
\begin{enumerate}
    \item \textbf{Feasibility stability:} no new feasibility cuts are generated during the backward pass;\vspace{-2pt}
    \item \textbf{Optimality:} the standard optimality gap criterion is satisfied, i.e., $|Z_{\mathrm{up}}-Z_{\mathrm{low}}|\le \varepsilon$.
\end{enumerate}
This guarantees that no additional feasibility information can alter the policy at the visited states and that the cost approximation has converged within the prescribed tolerance.

\subsection{Algorithm summary}

\vspace{1mm}
The overall procedure of the Penalty-Free SDDP algorithm is summarized in Algorithm~\ref{alg:pfsddp}. 
The pseudocode highlights the alternation between feasibility and optimality steps, as well as the convergence check based on the stability of feasibility cuts and the optimality gap criterion.

\begin{algorithm}[!htbp]
\caption{Penalty-Free SDDP}
\label{alg:pfsddp}
\begin{algorithmic}[1]
\State \textbf{Initialize:} empty sets of optimality and feasibility cuts
\Repeat
    \State $\mathcal{S}_t \gets \emptyset$, for all $t = 1{:}T$
    \State \textbf{Forward pass:}
    \For{each sampled path $\omega \in \Omega^{\mathrm{fwd}}$}
        \For{$t = 1$ to $T$}
            \State Solve~\eqref{eq:feasibility_problem_strict} at $(x_{t-1}(\omega),\xi_{t-1}(\omega))$; obtain $s_t^\ast(\omega)$ and $\beta_t^\ast$
            \State Solve~\eqref{eq:optimality_problem_with_FFF} with $0\le s_t \le s_t^\ast(\omega)$ and $\beta_t \le \beta_t^*$
            \State Record trial state: $\mathcal{S}_t \gets \mathcal{S}_t \cup \{x_t(\omega)\}$
        \EndFor
    \EndFor
    \State \textbf{Backward pass:}
    \For{$t = T$ down to $1$}
        \For{each $x_t \in \mathcal{S}_t$}
            \For{each $\xi_t \in \Xi_t$}
                \State Solve~\eqref{eq:feasibility_problem_strict} at $(x_t,\xi_t)$ to obtain $s_t^\ast$ and $\beta_t^\ast$ and add \textit{feasibility cuts} to update $\widetilde{\mathrm{FFF}}_{t}$
                \State Solve~\eqref{eq:optimality_problem_with_FFF} at $(x_t,\xi_t)$ with $0\le s_t \le s_t^\ast$ and $\beta_t \le \beta_t^\ast$ and add \textit{optimality cuts} to update $\widetilde{\mathrm{FCF}}_{t}$
            \EndFor
        \EndFor
    \EndFor
    \State Update $Z_{\mathrm{up}}$ and $Z_{\mathrm{low}}$
\Until{no new feasibility cuts added in backward pass \textbf{and} $|Z_{\mathrm{up}}-Z_{\mathrm{low}}|\le \varepsilon$}
\end{algorithmic}
\end{algorithm}

\vspace{-1mm}
This implementation ensures that the recursive process stop only when no additional feasibility information can modify the state decisions and the cost approximation has converged within tolerance.

\section{Case Study: Application to Energy Planning}

The proposed methodology was applied to a large-scale planning model inspired by the operational challenges of the Brazilian hydrothermal system. 
This system comprises dozens of interconnected reservoirs with complex physical interdependencies, where operators must reconcile multiple and often conflicting objectives — such as maximizing energy production, ensuring environmental and navigational flows, and meeting strict regulatory and ecological requirements. 
These \textit{hydraulic operating conditions} impose binding physical limits on the system and represent some of the most demanding aspects of large-scale hydrothermal scheduling, particularly when inflows are uncertain and their satisfaction competes with economic efficiency.

In current operational models, these conditions are typically represented either as hard constraints or through large artificial penalty terms. 
While the latter approach guarantees numerical feasibility, it introduces severe calibration challenges and distorts the interpretation of marginal costs. 
The Penalty-Free SDDP formulation proposed in this work was designed precisely to overcome these limitations, ensuring compliance with physical and environmental requirements without contaminating the economic objective with arbitrary penalties.

To provide a controlled validation of the proposed methodology, we focus here on a deterministic formulation of the problem. 
This setup can be solved as a single optimization problem that serves as a consistent reference solution, allowing the specific effects of feasibility cuts to be isolated from the stochastic variability that would otherwise arise in a fully stochastic implementation. 
Future developments will naturally extend this analysis to stochastic formulations, where inflow and renewable-generation uncertainties can further amplify the relevance of proper feasibility handling.

The test system represents a long-term hydrothermal scheduling model comprising 75 interconnected reservoirs arranged in cascades, 109 thermal plants, and 32 aggregated renewable plants. The planning horizon spans 4.5 years with monthly resolution, plus an additional buffer year to mitigate the so-called “End of the World” effect. 

Within the Brazilian hydrothermal context, these constraints define the physical, regulatory, and environmental operating conditions for hydropower plants. 
They include minimum total outflows, ecological flow requirements, navigability support and structural safety limits that ensure compliance with federal and state regulations. Altogether, they account for 662 \textit{hard constraints} that embody nonnegotiable regulatory and ecological priorities within system operation.

To evaluate the performance of the proposed approach, three formulations were compared:
\begin{enumerate}[(i)]
\item the traditional SDDP model using penalized slack variables to enforce HOCs;
\item a single large-scale deterministic problem serving as the benchmark solution; and
\item the proposed Penalty-Free SDDP model with feasibility cuts.
\end{enumerate}
All formulations were solved under identical convergence criteria of 0.5\%. The main computational results and performance comparisons are summarized in Table~\ref{table_performance}.

\begin{table}[!ht]
\renewcommand{\arraystretch}{1.3}
\centering
\caption{Performance comparison of the three formulations.}
\label{table_performance}
\begin{tabular}{|c|c|c|}
\hline
\textbf{Solution Method} & \textbf{Operation Cost [M\$]} & \textbf{Violation Cost [M\$]}\\
\hline
Single problem & 31,872.67 & 1.38 \\
\hline
SDDP & 31,877.09 & 2.59 \\
\hline
Penalty-free SDDP & 31,877.14 & 1.37 \\
\hline
\end{tabular}
\end{table}

\begin{figure}[!ht]
\centering
\includegraphics[width=\columnwidth]{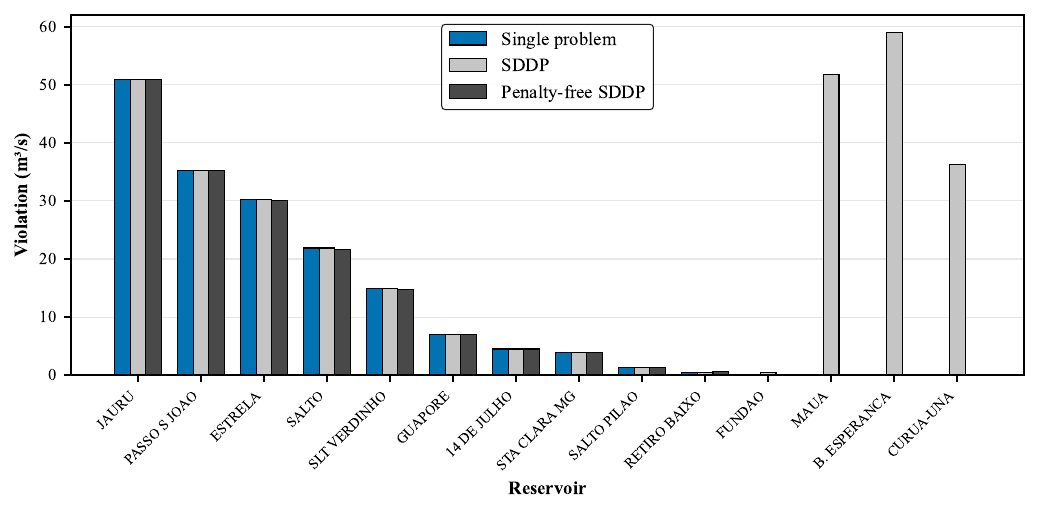}
\caption{Total minimum outflow violations.}
\label{fig_comparacao_violacoes_min_outflow}
\end{figure}

The Penalty-Free SDDP produced violations that were almost identical to those obtained in the single large-scale problem, affecting the same set of reservoirs, as shown in Figure~\ref{fig_comparacao_violacoes_min_outflow}. It correctly identified situations in which the strict enforcement of HOCs became infeasible under specific inflow conditions, reporting them explicitly through the feasibility function rather than masking them with large penalty terms. 
In contrast, the traditional SDDP showed violations at 14 reservoirs, compared to 11 in the single-problem benchmark, suggesting difficutly to calibrate the constraint and premature convergence. These results demonstrate that the use of feasibility cuts improves the interpretability and reliability of the solution, enabling the algorithm to distinguish between temporary and truly infeasible constraints without requiring iterative penalty tuning or repeated recalibration of model parameters.

Finally, regarding execution times, our proposed framework took 67 iterations to converge, requiring approximately 1.2$\times$ more iterations than the traditional SDDP to achieve convergence (56 iterations), translating to proportionally longer execution times. However, this comparison must consider that the SDDP produced some distortions in violation magnitudes and locations. Moreover, our framework directly produces implementable solutions with certified feasibility. For applications where satisfying hard system constraints is critical and understanding which constraints cannot be met provides essential information, the additional computational cost could be justified and the trade-off would be acceptable.

\section{Conclusions and Future Work}

This paper presented the \textit{Penalty-Free SDDP}, a general extension of the Stochastic Dual Dynamic Programming algorithm that eliminates the need for artificial penalty terms by incorporating feasibility cuts into its recursive structure. The proposed formulation separates economic optimization from constraint enforcement through an auxiliary feasibility recursion, which allows the algorithm to handle infeasibilities automatically and to explicitly identify \textit{truly infeasible} situations that cannot be corrected through adjustments in the state decisions.

The methodology was validated in a large-scale deterministic application inspired by the Brazilian hydrothermal system, where complex hydraulic, environmental, and regulatory conditions were modeled as hard constraints. The results demonstrated that the proposed approach preserves the economic consistency of the solution while correctly diagnosing infeasible operating conditions, without requiring iterative tuning or subjective calibration of penalty parameters — a process that is inherently case-dependent and often compromises model transparency.

A key advantage of this penalty-free scheme is its robustness and generality: by decoupling feasibility from cost minimization, it provides a principled way to ensure constraint satisfaction across different systems and formulations, while maintaining the interpretability of economic signals. 

Future research will focus on extending the method to fully stochastic, multi-stage settings — the natural domain of SDDP — where uncertainty in inflows and renewable generation amplifies the relevance of proper feasibility handling. 
Such an extension will inevitably increase the number of cuts, as feasibility cuts are generated for every stage, scenario and realization in the backward pass. 
Although this poses a computational challenge, it also removes the need for problem-specific penalty calibration and enables more accurate and transparent optimization processes in which feasibility and economic optimality are treated as distinct, complementary components of the decision framework.

\bibliographystyle{ieeetr}
\bibliography{bibliography.bib}






\end{document}